\documentclass[preprint,1p,times]{elsarticle}

\usepackage{amssymb}
\usepackage{amsmath}
\usepackage{amsfonts}
\usepackage{amssymb}
\usepackage{verbatim}
\usepackage{graphicx} 
\usepackage{color}

\newtheorem{theorem}{Theorem}

\newtheorem{definition}[theorem]{Definition}

\newtheorem{proposition}[theorem]{Proposition}
\newtheorem{remark}[theorem]{Remark}

\journal{Journal of Dinamical and Control Systems}

\begin{document}

\begin{frontmatter}

\title{The chain control set of discrete-time linear systems on the affine two-dimensional Lie group}

\author[inst1]{Thiago Matheus Cavalheiro}

\affiliation[inst1]{organization={Department of Mathematics, State University of Maringá},
            city={Maringá},
            state={Paraná},
            country={Brazil}}

\author[inst2]{João Augusto Navarro Cossich}
\author[inst3]{Alexandre José Santana}

\affiliation[inst2]{organization={Department of Mathematics, State University of Maringá},
            city={Maringá},
            state={Paraná},
            country={Brazil}}

\affiliation[inst3]{organization={Department of Mathematics, State University of Maringá},
            city={Maringá},
            state={Paraná},
            country={Brazil}}

\begin{abstract}

In this paper, we present  conditions for the existence and uniqueness of chain control sets of discrete-time linear systems on the affine two-dimensional Lie group. More specifically, we prove that these chain control sets are given by the union of an infinite number of control sets with  empty interiors.
\end{abstract}


\end{frontmatter}

\section{Introduction}
\label{sec:sample1}
A control system is controllable if, for all pairs of points or states in the state space, there is a trajectory of the system connecting them. If small jumps between pieces of trajectories are allowed to connect the points, then the system is chain controllable (see Colonius; Kliemann \cite{CCS2}).   However, it is well known that a system is not always (chain) controllable, and in several applications it is interesting to study subsets of the state space where this (chain) controllability property holds, known as (chain) control sets (see next section for precise definitions).
The theory related to controllability and chain controllability has been intensively studied over the past forty years. (e.g.  Albertini; Sontag \cite{sontalb} and \cite{sontalb1}, Ayala; Da Silva \cite{AyalaeAdriano}, Colonius; Cossich; Santana \cite{CCS} and \cite{CCS2} and Sontag \cite{sontag1}). 

In special, in \cite{CCS} the authors proved that for discrete-time linear systems on $\mathbb{R}^n$ with restricted control-range, the system admits a unique control set with no empty interior if and only if the unrestricted system is controllable. 
On the other hand, considering the state space as the two-dimensional affine group $\hbox{Aff}(2,\mathbb{R})$, the continuous case was completely studied in \cite{AyalaeAdriano}, that is,  all  linear systems on $\hbox{Aff}(2,\mathbb{R})$ were classified, the conditions for controllability were studied, and the authors also presented conditions for existence of control sets.
However, still in this state space, the global controllability problem in the discrete case was addressed by Cavalheiro; Cossich; Santana \cite{TAJ}, i.e., all possible discrete-time linear systems on $\hbox{Aff}(2,\mathbb{R})$ were described, and the conditions for controllability was studied.

In our work, we find all control sets of  discrete-time linear systems on $\hbox{Aff}(2,\mathbb{R})$ in  cases where the reachable sets have empty interiors. Roughly, the reachable set of a system is the set of all possible states that the system can reach from a given initial state over a specified time interval (see Section \ref{sec:sample:appendix}). Thereof, we present conditions for the existence and uniqueness of chain control sets. In fact, under some spectral conditions of the system we prove that the chain control set is the union of  an infinite number of control sets with empty interiors. Additionally, we study the chain controllability of the system.

This paper is divided as follows: Section 2 is dedicated to setting up notation, terminology, and some basic definitions. In sections 3 and 4, the precise expressions for the control sets are presented. Moreover, the existence and uniqueness of the chain control sets are proven in terms of the eigenvalues of the automorphism that defines the system. Additionally, in section 4, conditions for chain controllability of the system are studied.


\section{Preliminaries}
\label{sec:sample:appendix}
In this section, we establish the notations, definitions, and basic concepts that will be used throughout this paper. Consider a $n-$dimensional Riemannian manifold $M$ endowed with a metric $d$. Take a $\mathcal{C}^{\infty}$ map $f: U \times M \longrightarrow M$, where $U$ is a non-empty compact convex neighborhood of $0 \in \mathbb{R}^m$ such that $U \subset \overline{\hbox{int} U}$. Using the notation $f_u(x) = f(u,x)$, the system considered in this paper is defined by
\begin{equation}\label{sistgeneral}
    \Sigma: x_{k+1} = f_{u_k}(x_k), u_k \in U,
\end{equation}
with $k \in \mathbb{N}_0 = \mathbb{N} \cup \{0\}.$

Given an initial condition $x \in M$, the solution of $(\Sigma)$ is denoted by $\varphi(k,x,u)$ where $u \in \mathcal{U} = \prod_{i \in \mathbb{Z}} U$ such that $u = (u_i)_{i \in \mathbb{Z}}$. Assuming that $f_u$ is a diffeomorphism, the solution of $(\Sigma)$ is given by

\begin{equation*}
    \varphi(k,x,u) = 
    \left\{
    \begin{array}{ccl}
        f_{u_{k-1}} \circ \cdots \circ f_{u_0}(x),& k > 0.\\
        x,& k=0.\\
        f_{u_{k}}^{-1} \circ \cdots \circ f_{u_{-1}}^{-1}(x),& k < 0.
    \end{array}
    \right.
\end{equation*}
and, as is well-known in control theory (see \cite{CCS1}), the solution $\varphi$ satisfies the cocycle property. Therefore, the space $\mathcal{U}$ is compact. Now, define  $\Theta: \mathbb{Z} \times \mathcal{U} \longrightarrow \mathcal{U}$ by 
\begin{equation*}
    \Theta_k((u_j)_{j \in \mathbb{Z}}) = (u_{j+k})_{j \in \mathbb{Z}}. 
\end{equation*}

This map is continuous and therefore, $\Theta$ define a continuous dynamical system. Note that the solution $\varphi$ satisfies 
\begin{itemize}
    \item[1-] $\varphi(t + s, g, u) = \varphi(t, \varphi(s,g,u), \Theta_s(u)), t,s \in \mathbb{Z}$.
    \item[2-] If $ts > 0$ in $\mathbb{Z}$, $u,v \in \mathcal{U}$, then there is a $w \in \mathcal{U}$ such that 
    \begin{equation*}
        \varphi(t,\varphi(s,g,u),v) = \varphi(t + s, g, w), \forall g \in M.
    \end{equation*}
\end{itemize}

Recall that the set of points reachable and controllable from $x$ up to time $k \in \mathbb{N}$ are
\begin{eqnarray*}
    \mathcal{R}_k(x) &= \{y \in M: \hbox{ there is }u \in \mathcal{U} \hbox{ with }\varphi(k,x,u)=  y\}\\
    &\hbox{ and }\\ 
    \mathcal{C}_k(x) &= \{y \in M: \hbox{ there is }u \in \mathcal{U} \hbox{ with }\varphi(k,y,u)= x\}
\end{eqnarray*}

The sets $\mathcal{R}(x) = \bigcup_{k \in \mathbb{N}} \mathcal{R}_k(x)$ and $\mathcal{C}(x) = \bigcup_{k \in \mathbb{N}} \mathcal{C}_k(x)$ denote the reachable set and the controllable set from $x$ respectively. The next definition is one of the most important of this paper. 

\begin{definition}\label{defincontrolset}  For a system of the form (\ref{sistgeneral}), a nonvoid set $D \subset M$ is a control set if satisfies the following properties:
\begin{itemize}
    \item[1-] For all $x \in D$ there is a control $u \in \mathcal{U}$ such that $\varphi(k,x,u) \in D$ for all $k \in \mathbb{N}$,
    \item[2-] For any $x \in D$ one has $D \subset \overline{\mathcal{R}(x)}$,
    \item[3-] $D$ is maximal with such properties.  
\end{itemize}
\end{definition}

The maximality property guarantees that if $D' \supset D$  satisfies conditions $1$ and $2$ of the above definition, then $D' = D$. Besides, it is not hard to prove that if two control sets $D$ and $D'$ satisfies $D \cap D' \neq \emptyset$, then $D = D'$. The second property is usually called approximate controllability. Such property can be weakened, if we consider only "pieces of trajectory", as in the next definition.

\begin{definition}\label{controlablechain}Given $x,y \in M$, $k \in \mathbb{N}$ and a real number $\varepsilon>0$, an $(\varepsilon,k)-$controllable chain between $x$ and $y$ is given by the set 
\begin{equation*}
    \mu_{(\varepsilon,k)} = 
    \begin{Bmatrix}
    n \in \mathbb{N}\\
    x_0,x_1, \ldots ,x_n \in M\\
    u_0, \ldots ,u_{n-1} \in \mathcal{U}\\
    k_0, \ldots ,k_{n-1} \in \mathbb{N} \setminus \{1, \ldots ,k-1\}
    \end{Bmatrix}
\end{equation*}
satisfying $x_0 = x$, $y = x_n$ and 
 \begin{equation*}
    d(\varphi(k_i,x_i,u_i),x_{i+1}) < \varepsilon, \forall i=0, \ldots ,n-1. 
\end{equation*}

We say that the pair $x, y \in M$ is chain controllable if for all  $(\varepsilon,k) \in \mathbb{R}^+ \times \mathbb{N}$, there is a $(\varepsilon,k)-$controllable chain between $x$ and $y$. 
\end{definition}

When we consider the maximal subset of $M$ with chain controllability property, we have the following definition. 

\begin{definition}\label{definchaincontrolset} A nonvoid set $E \subset M$ is a chain control set of (\ref{sistgeneral}) if satisfies the following properties:
\begin{itemize}
    \item[1-] For all $x \in E$, there is a control $u \in \mathcal{U}$ such that $\varphi(k,x,u) \in E$, for all $k \in \mathbb{Z}$,
    \item[2-] Every pair $x,y \in E$ is chain controllable.   
    \item[3-] $E$ is maximal with such properties.  
\end{itemize}
\end{definition}

From now on, let us consider $M = G$ a connected $n-$dimensional Lie group. Then, in the next definition we specify the concept of discrete linear control system (see e.g. \cite{CCS1}) for Lie groups.   

\begin{definition}\label{lindefin} Consider the discrete-time control system on $G$
\begin{equation*}
    \Sigma: x_{k+1} = f_{u_k}(x_k), u_k \in U,
\end{equation*}
where $U \subset \mathbb{R}^m$ is a compact convex neighborhood of $0$. Then ($\Sigma$) is called linear if $f_0: G \longrightarrow G$ is an automorphism and 
\begin{equation}\label{transl}
    f_u(g) = f_u(e) \cdot f_0(g). 
\end{equation}
where $"\cdot"$ denotes the operation of $G$.
\end{definition}

The map $f$ can be defined using the translations. Fixing $u \in U$ and knowing that $f_u(e) \in G$, the expression (\ref{transl}) allow us to write $f_u(g)$ as 
\begin{equation}\label{transl1}
    f_u(g) = f_u(e)f_0(g) = L_{f_u(e)}(f_0(g)), 
\end{equation}
where $L_{f_u(e)}$ is the left translation by the element $f_u(e)$. Considering the expression above, the inverse of $f_u$ is given by 
\begin{equation}
    (f_u)^{-1}(g) = f_0^{-1}\circ L_{(f_u(e))^{-1}}(g)= f_0^{-1}((f_u(e))^{-1} \cdot g). 
\end{equation}

Then, we can conclude that $f_u$ is a diffeomorphism of $G$, for any $u \in U$. The next result shows that the solutions can also be defined in terms of translations (see \cite[Proposition 3]{CCS1}). 

\begin{proposition}\label{prop52} Consider a discrete-time linear control system on a Lie group $G$ defined by $x_{k+1} = f(u_k,x_k)$, $u_k \in U$. Then for all $g \in G$ and $u = (u_i)_{i \in \mathbb{Z}} \in \mathcal{U}$ we have
\begin{equation*}
    \varphi(k,g,u) = \varphi(k,e,u)f_0^k(g). 
\end{equation*}
\end{proposition}

\section{Control sets of $\Sigma$} 

It is well known that the affine two-dimensional Lie group given by 
\begin{equation*}
    \hbox{Aff}(2,\mathbb{R})= \left\{\begin{bmatrix}
        x & y \\
        0 & 1
    \end{bmatrix}: x > 0, y \in \mathbb{R}\right\},
\end{equation*}
can be considered as the half-plane $G = \mathbb{R}^+ \times \mathbb{R}$ endowed with the semidirect product 
\begin{equation*}
    (x,y) \cdot (z,w) = (xz, w + yz). 
\end{equation*}

Following the definition (\ref{lindefin}) and according to \cite{AyalaeAdriano}, every automorphism of $\hbox{Aff}(2,\mathbb{R})$ has the form 
\begin{equation*}
    \phi(x,y) = (x,a(x-1) + dy), \,\, d \neq 0. 
\end{equation*}

From this we have that the discrete-time linear control systems on $\hbox{Aff}(2,\mathbb{R})$ are given by 
\begin{equation*}
    \Sigma: (x_{k+1},y_{k+1}) = f_{u_k}(x_k,y_k), k \in \mathbb{N},
\end{equation*}
with  $f: U \times \hbox{Aff}(2,\mathbb{R}) \longrightarrow \hbox{Aff}(2,\mathbb{R})$  defined by 
\begin{equation*}
    f_u(x,y) = (h(u)x, a(x-1) + dy + g(u)x), \,\, d \neq 0
\end{equation*}
where the maps $h: U \longrightarrow \mathbb{R}^+$ and $g:U \longrightarrow \mathbb{R}$ are $\mathcal{C}^{\infty}$ and satisfy $h(0)=1$ and $g(0)=0$ with $U$ being  a compact convex neighborhood of $0 \in \mathbb{R}^m$. 

Noting that 
\begin{equation*}
    d(f_0)_{(1,0)} = 
    \begin{pmatrix}
    1 & 0 \\
    a & d
    \end{pmatrix}
\end{equation*}
we have  $\hbox{Spec}\{(df_0)_{\beta}\} = \{1,d\}$ for any basis $\beta$ of $\mathbb{R}^2$. In the current case, by taking $h(u) \equiv 1$, the solution has the form 
\begin{equation*}
    \varphi(k,(x,y),u) = (x, S^k_u(x,y)), 
\end{equation*}
where $S_k^u(x,y) = \sum_{j=0}^{k-1} d^{k-1-j}[(a+g(u_j))x - a] + d^k y$.

The goal of this section is to explicitly describe the control sets of the above system. One can easily check that for any $(x,y) \in \hbox{Aff}(2,\mathbb{R})$ we have $\varphi(k,(x,y),u) \in \{x\} \times \mathbb{R}$ for every $u \in \mathcal{U}$ and $k \in \mathbb{Z}$. Then, $\hbox{int}\mathcal{R}(x,y) = \emptyset,$ for every $(x,y) \in \hbox{Aff}(2,\mathbb{R})$. In particular, this implies the non-controllability of the system $(\Sigma)$. Let us describe the control sets of $(\Sigma)$. 

To obtain  the control sets we consider the eigenvalue $d$ such that $|d| < 1$. First, consider $0 < d < 1$. Note that $\{[(a+g(u))x - a]: u \in U\}$ is a compact connected subset of $\mathbb{R}$. In particular, there are $u_m, u_M \in \mathcal{U}$ such that 
\begin{eqnarray*}
    M(x) &=& \max_{u \in U}\{(a+g(u))x - a\} = (a+g(u_M))x - a.\\
    m(x) &=& \min_{u \in U}\{(a+g(u))x - a\} = (a+g(u_m))x - a.
\end{eqnarray*}

Hence the description of the control sets in case of $0 < d < 1$ can be summarized in the following proposition.

\begin{proposition}\label{controlsetdpos}For any $x \in \mathbb{R}^+$ and $0<d<1$, the sets 
\begin{equation}\label{form1}
    D_x = \{x\} \times \left[\frac{m(x)}{1-d}, \frac{M(x)}{1-d}\right]. 
\end{equation}
are control sets of $(\Sigma)$. 
\end{proposition}

\noindent\textit{Proof:} Given $y \in \left[\frac{m(x)}{1-d}, \frac{M(x)}{1-d}\right]$, for any $u \in \mathcal{U}$ we get 
\begin{eqnarray*}
    S_k^u(x,y) &=& \sum_{j=0}^{k-1} d^{k-1-j} [(a+ g(u_j))x -a] + d^k y \leq \sum_{j=0}^{k-1} d^{k-1-j} M(x) + d^k y \\
    &=& \frac{1-d^k}{1-d}M(x) + d^k y \leq  \frac{1-d^k}{1-d}M(x) + d^k \left(\frac{M(x)}{1-d}\right)\\
    &=&\frac{M(x)}{1-d}. 
\end{eqnarray*}

Similarly, we get $S_k^u(x,y) \geq \frac{m(x)}{1-d}$. Then $S_k^u(x,y) \in \left[\frac{m(x)}{1-d}, \frac{M(x)}{1-d}\right]$ for any $k \in \mathbb{N}$ and $u \in \mathcal{U}$. Hence 
\begin{equation*}
    \varphi(k,(x,y),u) = (x, S_k^u(x,y)) \in \{x\} \times \left[\frac{m(x)}{1-d}, \frac{M(x)}{1-d}\right], \forall (k,u) \in \mathbb{N} \times \mathcal{U}.
\end{equation*}

To prove the approximate controllability, take $z,y \in \left[\frac{m(x)}{1-d}, \frac{M(x)}{1-d}\right]$. Then there is $u_y \in \mathcal{U}$ such that 
\begin{equation*}
    y = \frac{1}{1-d}[(a+g(u_y))x - a]. 
\end{equation*}

Taking $u = (u_y)_{j \in \mathbb{N}_0}$, then for all $\varepsilon>0$, there is a $k \in \mathbb{N}$ such that 
\begin{equation*}
    |d^k z | < \frac{\varepsilon}{2} \hbox{ and } \left|- \frac{d^k((a+g(u_y))x-a)}{1-d} \right| < \frac{\varepsilon}{2}. 
\end{equation*}

For $k$ and $u$ as before, we get
\begin{eqnarray*}
    || \varphi(k,(x,z),u) - (x,y)||
    &=&|\sum_{j=0}^{k-1}d^{k-1-j} [(a+g(u_y))x- a] + d^k z - y|\\
    &=&\left|\frac{1-d^k}{1-d}((a+g(u_y))x-a) + d^k z -\frac{1}{1-d}[(a+g(u_y))x - a]\right|\\
    &\leq&|d^k z| + \left|- \frac{d^k((a+g(u_y))x-a)}{1-d} \right| < \varepsilon. 
\end{eqnarray*}
which means $(x,z) \in \overline{\mathcal{R}(x,y)}. $

To prove the maximality, suppose $z > \frac{M(x)}{1-d}$, for each $(x,y) \in D_x$. As $y \in [S_k^{v_m}(x,y),S_k^{v_M}(x,y)] \leq \frac{M(x)}{1-d}$ and $S^u_k(x,y) \leq S_k^{v_M}(x,y) < z$, we choose $\varepsilon = z - \frac{M(x)}{1-d}$, to get 
\begin{equation*}
    S^u_k(x,y) \notin (z- \varepsilon, z+\varepsilon) = \left(\frac{M(x)}{1-d}, 2z - \frac{M(x)}{1-d}\right). 
\end{equation*}

Then $(x,z) \notin \overline{\mathcal{R}(x,y)}$ for any pair $(x,y) \in D_x$. The same argument can be used for $z < \frac{m(x)}{1-d}$. Hence $D_x$ is a control set for the system $(\Sigma)$. 
$\blacksquare$

\begin{figure}
    \centering
    \includegraphics[scale=0.6]{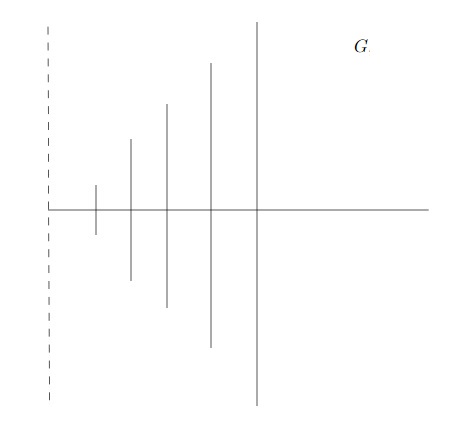}
    \caption{A representation of the control sets in case of $g(u) = u$, $a=0$ and $U = [-1,1]$. In the general case, the point $(x,0)$ does not need necessarily to be in the relative interior  (by the induced topology) of $D_x$.}
    \label{controlsets0<d<1}
\end{figure}

\begin{remark}With the same hypothesis above, supposing $D \subset G$ a control set of $(\Sigma)$ and $(x,y) \in D$, we have $D \subset \{x\} \times \mathbb{R}$, since $\varphi(k,(x,y),u)=(x,S_k^u(x,y))$. As $D_x$ is the maximal region of approximate controllability, then $D = D_x$. Therefore, if $D \subset G$ is a control set of $(\Sigma)$, then $D = D_x$, for some $x \in \mathbb{R}^+$. 
\end{remark}

Considering the map $p: U  \times \hbox{Aff}(2,\mathbb{R}) \longrightarrow \mathbb{R}$ defined by $p(u,x) = [(a+g(u))x-a]$, and supposing $-1 < d < 0$, it follows that  $\sum_{j=0}^{k-1} d^{k-1-j}p(u_j,x)$ satisfies
\begin{equation}\label{eq1}
    dM(x)\left(\frac{1-d^{k-1}}{1-d^2}\right) +m(x)\left(\frac{1-d^{k+1}}{1-d^2}\right) \leq \sum_{j=0}^{k-1} d^{k-1-j}p(u_j,x), 
\end{equation}
and 
\begin{equation}\label{eq2}
        \sum_{j=0}^{k-1} d^{k-1-j}p(u_j,x) \leq  M(x)\left(\frac{1-d^{k+1}}{1-d^2}\right) +dm(x)\left(\frac{1-d^{k-1}}{1-d^2}\right)
\end{equation}
if $k \geq 3$ is odd. If $k \geq 2$ is even, we have 
\begin{equation}\label{eq3}
    (dM(x) +m(x))\frac{1-d^k}{1-d^2}\leq \sum_{j=0}^{k-1}d^{k-1-j} p(u_j,x) \leq (M(x) +dm(x))\frac{1-d^k}{1-d^2}.
\end{equation}

Therefore, the control sets in  case of $-1 < d < 0$ are characterized by the next proposition. 

\begin{proposition} If $-1 < d < 0$ and $h \equiv 1$, for each $x \in \mathbb{R}^+$, then the sets 
\begin{equation}\label{form2}
    D_x = \{x\} \times \left[\frac{d M(x) + m(x)}{1-d^2}, \frac{M(x) + dm(x)}{1-d^2}\right]. 
\end{equation}
are control sets of the system $(\Sigma)$. 
\end{proposition}

\noindent\textit{Proof:} Take $x \in \mathbb{R}^+$ and $y \in \left[\frac{d M(x) + m(x)}{1-d^2}, \frac{M(x) + dm(x)}{1-d^2}\right]$. For $k$ even, we have
\begin{eqnarray*}
    S_k^u(x,y) &=& \sum_{j=0}^{k-1} d^{k-1-j} p(u_j,x) + d^k y
    \leq  (M(x) +dm(x))\frac{1-d^k}{1-d^2} + d^k y\\
    &\leq&  (M(x) +dm(x))\frac{1-d^k}{1-d^2} + d^k \left(\frac{M(x) + dm(x)}{1-d^2}\right) =\frac{M(x) + dm(x)}{1-d^2} 
\end{eqnarray*}
for all $k\geq 2$ even and $u \in \mathcal{U}$. Similarly we prove
\begin{equation*}
     S_k^u(x,y) \geq \frac{d M(x) + m(x)}{1-d^2}.
\end{equation*}

For $k= 1$, we have
\begin{equation*}
    S_1^{u}(x,y) = p(u_0,x) + dy \leq M(x) + dy. 
\end{equation*}

As $d < 0$ it follows that 
\begin{equation}\label{d_k_odd}
    d^k\left(\frac{M(x) + dm(x)}{1-d^2}\right)\leq d^k y \leq d^k\left(\frac{dM(x) + m(x)}{1-d^2}\right).
\end{equation}
for every $k \in \mathbb{N}$ odd. Then
\begin{eqnarray*}
    M(x) + dy &\leq& M(x) + d\left(\frac{dM(x) + m(x)}{1-d^2}\right)
    = M(x) \left(\frac{1-d^2}{1-d^2}\right) +\\
    & &+\left(\frac{d^2M(x) + dm(x)}{1-d^2}\right)\\
    &=& \frac{M(x) + dm(x)}{1-d^2} 
\end{eqnarray*}
implying  that
\begin{equation*}
    S_1^u(x,y) \leq \frac{M(x) + dm(x)}{1-d^2}
\end{equation*}
for all $u \in \mathcal{U}$. By a  similar argument, we obtain 
\begin{equation*}
    S_1^u(x,y) \geq \frac{dM(x) + m(x)}{1-d^2}
\end{equation*}
for any $u \in \mathcal{U}$. 

Now, suppose $k \geq 3$ odd. By (\ref{d_k_odd}) we get 
\begin{eqnarray*}
    S_k^u(x,y) &=& \sum_{j=0}^{k-1} d^{k-1-j} p(u_j,x) + d^k y 
    \leq \frac{M(x) + dm(x)}{1-d^2} -\\
    & &- d^k \left(\frac{dM(x) + m(x)}{1-d^2}\right) + d^k y
    \leq \frac{M(x) + dm(x)}{1-d^2}. 
\end{eqnarray*}

Likewise, 
\begin{equation*}
    S_k^u(x,y) \geq \frac{dM(x) + m(x)}{1-d^2} 
\end{equation*}
for every $k \geq 3$ odd and $u \in \mathcal{U}$. Thus $\varphi(k,(x,y),u) \in D_x$, for every $k \in \mathbb{N}$ and $u \in \mathcal{U}$. Now we show the approximate controllability. Given $x \in \mathbb{R}^+$ and $y,z \in \left[\frac{d M(x) + m(x)}{1-d^2}, \frac{M(x) + dm(x)}{1-d^2}\right]$, there are $u,v \in U$ such that 
\begin{equation*}
    z = \frac{dp(u,x) + p(v,x)}{1-d^2}. 
\end{equation*}

Taking $\varepsilon>0$, there is a $k \in \mathbb{N}$ even such that 
\begin{equation*}
    |d^k y | < \frac{\varepsilon}{2} 
\end{equation*}
and 
\begin{equation*}
    \left|-\frac{d^k(dp(u,x) + p(v,x))}{1-d^2} \right| < \frac{\varepsilon}{2}. 
\end{equation*}

Then, taking $w_{2j} = u$ and $w_{2j+1} = v$ for all $j \in \mathbb{N}_0$ and $w = (w_j)_{j \in \mathbb{N}_0}$ we get
\begin{eqnarray*}
    ||\varphi(k,(x,y),w) - (x,z)|| &=& |\sum_{j=0}^{k-1} d^{k-1-j}p(w_j,x) + d^k y - z|\\
    &=&\left|\sum_{j=0}^{\frac{k}{2}-1} d^{k-1-2j} p(w_{2j},x) + \sum_{j=0}^{\frac{k}{2}-1} d^{k-2-2j} p(w_{2j+1},x) + d^k y - z  \right|\\
    &=&\left| d\left(\frac{1-d^k}{1-d^2}\right) p(u,x) + \left(\frac{1-d^k}{1-d^2}\right) p(v,x) + d^k y - z\right|\\
    &=&\left|  \frac{dp(u,x) + p(v,x)}{1-d^2} - d^k\left( \frac{dp(u,x) + p(v,x)}{1-d^2}\right) + d^k y- z\right|\\
    &<& \varepsilon. 
\end{eqnarray*}

Then $(x,y) \in \overline{\mathcal{R}(x,z)}$. The maximality follows using the same argument as in Proposition (\ref{controlsetdpos}). $\blacksquare$

The next section is focused on the chain control sets associated with the system $(\Sigma)$. 

\section{Chain control sets of $\Sigma$} 

As the previous section, we consider the   eigenvalue $d$ such that $|d| < 1$. Take the maps $\eta, \,\, \mu: \mathbb{R}_+ \longrightarrow \mathbb{R}$ given by 
\begin{equation*}
    \eta(x) = \frac{a +g(u_M)}{1-d} x - \frac{a}{1-d} \hbox{ and }\mu(x) = \frac{a +g(u_m)}{1-d} x - \frac{a}{1-d}.
\end{equation*}

Note that the region between their traces defines a surface in $G$, given by 
\begin{equation*}
    E = \left\{(z,t) \in G:\frac{a +g(u_m)}{1-d} z - \frac{a}{1-d} \leq t \leq \frac{a +g(u_M)}{1-d} z - \frac{a}{1-d} \right\}.
\end{equation*}

The next step is to prove that $E = \bigcup_{x \in \mathbb{R}^+} D_x$, where $D_x$ depends on the sign of $d$ and can be defined by (\ref{form1}) (case $0 < d < 1$) or (\ref{form2}) (case $-1 < d < 0$). Let us start with the case $0 < d < 1$. 

\begin{proposition}\label{chaincontrolsetd>0} Considering the case $0 < d < 1$ and $h \equiv 1$, the set
\begin{equation*}
    E = \bigcup_{x \in \mathbb{R}_+} D_x, 
\end{equation*}
with $D_x$ defined in the equality (\ref{controlsetdpos}), is the unique chain control set of $(\Sigma)$. 
\end{proposition}

\noindent\textit{Proof:} Note that 
\begin{equation*}
    f^{-1}_u(x,y) = \left(x, \frac{-a}{d}\left(x - 1\right) + \frac{1}{d}\left(y - g(u)x\right)\right).
\end{equation*}

Using the fact that there is an $\alpha \in U$ such that $y = \frac{p(\alpha,x)}{1-d}$ and taking $u = (\alpha)_{i \in \mathbb{Z}}$, we get
\begin{eqnarray*}
    \frac{1}{d}(-a(x-1) +y - g(\alpha)x) = \frac{1}{d}\left(y - p(\alpha,x)\right) = \frac{1}{d}\left(\frac{p(\alpha,x)}{1-d} - p(\alpha,x)\right)=\frac{1}{d}\left( \frac{d p(\alpha,x)}{1-d}\right)=\frac{p(\alpha,x)}{1-d} = y,
\end{eqnarray*}
which means that $\varphi(-k,(x,y),u) = (x,y) \in \{x\}\times \left[\frac{m(x)}{1-d},\frac{M(x)}{1-d}\right]$, for any $k \in \mathbb{Z}$.       

To prove the chain controllability of $E$, take $p=(x,y)$ and $q = (z,t)$ in $E$. If $x=z$, one can use the approximate controllability to show the chain controllability between $p$ and $q$. Without loss of generality, suppose that $x < z$ and take $(\varepsilon, k) \in \mathbb{R}_+ \times \mathbb{N}$. There are $T_0 \geq k$ in $\mathbb{N}$ and $u_0 \in \mathcal{U}$ such that $\varphi(T_0,(x,y), u_0) \in B_{\varepsilon}(x,\frac{M(x)}{1-d})$. If $z \in (x, x + \varepsilon)$, we can take $(z,y') \in \{z\} \times \mathbb{R}$ such that $\varphi(T_0,(x,y), u_0) \in B_{\varepsilon}(z,y')$. If $(z,y') \notin D_z$, we take $u_1 = (u_M)_{i \in \mathbb{Z}}$ and $T_1 \geq k$ such that $d^{T_1} y < \frac{\varepsilon}{2}$ and $-d^{T_1}\left(\frac{M(z)}{1-d}\right) < \frac{\varepsilon}{2}$, then $S_{u_1}^{T_1}(x,y) < \frac{M(z)}{1-d}+\varepsilon$ and hence  
\begin{equation*}
    \varphi(T_1,(z,y'),u_1) \in B_{\varepsilon
    }\left(z,\frac{M(z)}{1-d}\right).
\end{equation*}

Then, by approximate controllability of $D_z$, there are $T_2 \geq k$ and $u_2 \in \mathcal{U}$ such that 
\begin{equation*}
    \varphi(T_2, (z,\frac{M(z)}{1-d}),u_2) \in B_{\varepsilon}(z,t),
\end{equation*}
which proves the chain controllability between $p$ and $q$. If $z \geq x+\varepsilon$, there is a $x_1 \in (x,x+\varepsilon)$ and $y_1 \in \mathbb{R}$ such that $\varphi(T_0,(x,y),u_0) \in B_{\varepsilon}(x_1,y_1)$. Taking $u_1 = (u_M)_{i \in \mathbb{Z}}$ and $T_1 \geq k$ satisfying $d^{T_1}y_1 < \varepsilon$ and $-d^{T_1}\left(\frac{M(x_1)}{1-d}\right) < \frac{\varepsilon}{2}$ we get $\varphi(T_1,(x_1,y_1),u_1) \in B_{\varepsilon}(x_1,\frac{M(x_1)}{1-d}).$ If $z \in (x_1,x_1 + \varepsilon),$  take $(z,y') \in \{z\} \times \mathbb{R}$, $T_2 \geq k$ and $u_2 = (u_M)_{i \in \mathbb{Z}}$ satisfying $\varphi(T_2, (z,y'),u_2) \in B_{\varepsilon}(z,\frac{M(z)}{1-d})$. Using again the approximate controllability of $D_z$, there are $T_3 \geq k$ and $u_3 \in \mathcal{U}$ such that $\varphi(T_3,(z,\frac{M(z)}{1-d}),u_3) \in B_{\varepsilon}(z,t)$. Then, there exists $(x_{n-1},y_{n-1}) \in G$ such that $z \in (x_{n-1},x_{n-1} + \varepsilon)$, $T_{n-1} \geq k$ and $u_{n-1} \in \mathcal{U}$ satisfies $\varphi(T_{n-1},(x_{n-1},y_{n-1}),u_{n-1}) \in B_{\varepsilon}(z,t)$. Then, the set 
\begin{equation*}
    \mu_{(\varepsilon,k)} = 
    \left\{
    \begin{array}{ccc}
        (x_0,y_0) = (x,y), (x_1,y_1),...,(x_n,y_n) = (z,t) \in G\\
        u_0, u_1,..., u_{n-1} \in \mathcal{U} \\
        T_0,T_1,...,T_{n-1} \in \mathbb{N}\setminus \{0,...,k-1\}
    \end{array}    
    \right.
\end{equation*}
is a controlled chain between $p = (x,y)$ and $q = (z,t)$. 

Now, to prove the maximality of $E$, we show that if $(x,y) \notin E$, then there is a pair $(\delta,k) \in \mathbb{R}_+ \times \mathbb{N}$ such that no controlled chain with initial point in $E$ can reach $(x,y)$. Without loss of generality, consider the case $(x,y) \notin E$ with $y > \frac{M(x)}{1-d}$. Define $r = d((x,y),r_0)$, where
\begin{equation*}
    r_0 = \left\{(m,n) \in G: n = \frac{M(m)}{1-d}\right\}. 
\end{equation*}

 Taking $\delta = \frac{r}{4}$ and $0 < \varepsilon_0 < \delta,$ there is a $k \in \mathbb{N}$ such that $d^T(\delta + \varepsilon_0) < \varepsilon_0$ for any $T \geq k$. Now consider the region of $G$ defined by 
\begin{equation*}
    A_1 = \left\{(m,n) \in G: \frac{M(m        )}{1-d} < n \leq \frac{M(m)}{1-d} + \frac{\delta + \varepsilon}{\cos{\alpha}}\right\},
\end{equation*}
where $\alpha$ is the angle between the line containing $(x,y)$ orthogonal to $r_0$ and a line orthogonal to $\mathbb{R}_+ \times \{0\}$. One can prove that $\alpha \in [0,\frac{\pi}{2})$. Suppose by contradiction that the set 
\begin{equation*}
    \mu_{(\delta,k)} = 
    \begin{Bmatrix}
        p_0, p_1,...,p_n = (x,y) \in G\\
        u_0, u_1, ..., u_{n-1} \in \mathcal{U} \\
        T_0,T_1,..., T_{n-1} \in \mathbb{N} \setminus \{0,...,k-1\}
    \end{Bmatrix},
\end{equation*}
is a controlled chain between $p_0 \in E$ and $p_n=(x,y)$. Let $j_0 \in \{1,...,n-1\}$ be the smallest index such that $p_{j_0+i} \notin E$ and $p_{j_0} \in A_1$, with $i \in \{0,...,n-j_0\}$. Such index exists since the last point of the controllable chain $\mu_{(\delta,k)}$ is in the region $A_1$. Then $\varphi(T_{j_0-1},p_{j_0-1},u_{j_0-1}) \in B_{\delta}(p_{j_0})$, with $p_{j_0} \notin E$ and $p_{j_0-1} \in E$.  Considering $p_{j_0}= (p_{j_0}^1,p_{j_0}^2)$, we have 
\begin{eqnarray*}
    \sum_{j=0}^{T_{j_0}-1} d^{T_{j_0} - 1 - j} p(u_j,p_{j_0}^1) + d^{T_{j_0}} p_{j_0}^2 &\leq& \frac{1-d^{T_{j_0}}}{1-d}M(p_{j_0}) + d^{T_{j_0}} \left(\frac{M(p_{j_0})}{1-d} + \frac{\delta + r}{\cos{\alpha}}\right)\\
    &=&\frac{M(p_{j_0})}{1-d} + d^{T_{j_0}} \left(\frac{\delta + r}{\cos{\alpha}}\right) < \frac{M(p_{j_0})}{1-d} + \frac{\varepsilon}{\cos{\alpha}} , 
\end{eqnarray*}
that is
\begin{equation*}
    \varphi(T_{j_0},p_{j_0},u_{j_0}) \in B_{\delta}(p_{j_0+1}),
\end{equation*}
with $p_{j_0+1} \in A_1$. Similarly, we get  
\begin{eqnarray*}
    \sum_{j=0}^{T_{j_0+1}-1} d^{T_{j_0+1} - 1 - j} p(u_j,p_{j_0+1}^1) + d^{T_{j_0+1}} p_{j_0+1}^2 &\leq& \frac{1-d^{T_{j_0+1}}}{1-d}M(p_{j_0+1}) + d^{T_{j_0+1}} \left(\frac{M(p_{j_0+1})}{1-d} + \frac{\delta + r}{\cos{\alpha}}\right)\\
    &=&\frac{M(p_{j_0+1})}{1-d} + d^{T_{j_0+1}} \left(\frac{\delta + r}{\cos{\alpha}}\right) < \frac{M(p_{j_0+1})}{1-d} + \frac{\varepsilon}{\cos{\alpha}}. 
\end{eqnarray*}

Then $\varphi(T_{j_0+1},p_{j_0+1},u_{j_0+1}) \in B_{\delta}(p_{j_0+2})$, with $p_{j_0+2} \in A_1$. Such property is satisfied for any $p_{j_0+i}$, with $i \in \{0,...,n-j_0-1\}$, therefore, $\varphi(T_{n-1},p_{n-1},u_{n-1}) \notin B_{\delta}(x,y)$, contradiction. The same holds in case $(x,y) \notin E$, with $n < \frac{m(x)}{1-d}$. Thus, the set $E$ is the maximal region of chain controllability. The uniqueness is obvious. $\blacksquare$

\begin{figure}
    \centering
    \includegraphics[scale=0.6]{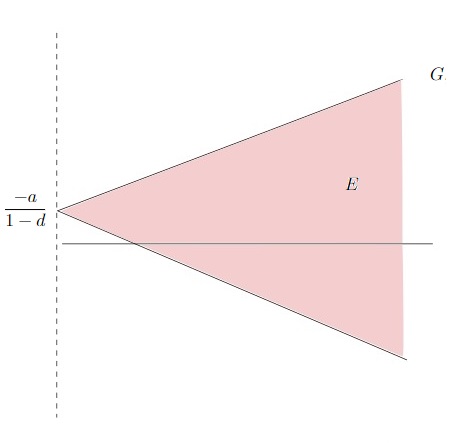}
    \caption{A representation of the chain control set of the Proposition \ref{chaincontrolsetd>0}, when $a < 0$.}
    \label{chaincontrolsetsd>0}
\end{figure}

\begin{proposition} Considering the case $-1 < d < 0$ and $h \equiv 1$, the set
\begin{equation*}
    E = \bigcup_{x \in \mathbb{R}_+} D_x, 
\end{equation*}
with $D_x$ defined in (\ref{form2}), is the unique chain control set of $(\Sigma)$. 
\end{proposition} 

\noindent\textit{Proof:} Consider $(x,y) \in D_x.$ Then, there are $u_1, u_2 \in U$ such that 
\begin{equation*}
    y= \frac{dp(u_1,x) + p(u_2,x)}{1-d^2}.
\end{equation*}

Taking $u_2 \in U,$ we have
\begin{eqnarray*}
    \frac{1}{d}\left(y - p(u_1,x)\right) &=&  \frac{1}{d}\left(\frac{dp(u_1,x) + p(u_2,x)}{1-d^2} - p(u_1,x)\right)\\
    &=&\frac{1}{d}\left(\frac{dp(u_1,x) + p(u_2,x)}{1-d^2} - p(u_1,x)\frac{1-d^2}{1-d^2}\right)\\
    &=&\frac{p(u_1,x) + dp(u_2,x)}{1-d^2}
\end{eqnarray*}
which, in particular satisfies $\dfrac{p(u_1,x) + dp(u_2,x)}{1-d^2} \in \left[\dfrac{m(x) + dM(x)}{1-d^2},\dfrac{M(x) + dm(x)}{1-d^2}\right]$. Now, taking $u_1 \in U$, we have
\begin{eqnarray*}
    \frac{1}{d}\left(\frac{p(u_1,x) + dp(u_2,x)}{1-d^2} - p(u_1,x)\right) &=&  \frac{1}{d}\left(\frac{p(u_1,x) + dp(u_2,x)}{1-d^2} - p(u_1,x)\frac{1-d^2}{1-d^2}\right)\\
    &=&\frac{1}{d}\left(\frac{dp(u_2,x)+d^2p(u_1,x)}{1-d^2}\right)=\frac{p(u_2,x) + dp(u_1,x)}{1-d^2}=y.
\end{eqnarray*}

We conclude that for $u = (...,u_2,u_1,u_2,...) \in \mathcal{U}$, we have $\varphi(-k,(x,y),u) \in D_x$ for all $k \in \mathbb{N}$. The chain controllability can be proved similarly to that in the previous result. For the maximality, take $(x,y) \in G$ such that $y > \frac{M(x) + dm(x)}{1-d^2}$. Define $r = d((x,y),r_0)$, with 
\begin{equation*}
    r_0 = \left\{(z,t) \in G: t = \frac{M(z) + dm(z)}{1-d^2}\right\}. 
\end{equation*}

Also consider $\delta = \frac{r}{4}$, $0 < \varepsilon < \delta$ and $\alpha \in [0,\frac{\pi}{2})$, where $\alpha$ is the angle between  $r_0$ and a line orthogonal to the axis $Ox$. Take $k \in \mathbb{N}$ such that for any $T \geq k$ we have $d^T(\varepsilon + \delta) < \varepsilon$. We claim that no $(\delta, k)-$controlled chain starting in $E$ can reach $(x,y)$. In fact, suppose by contradiction that the following set 
\begin{equation*}
    \mu_{(\delta,k)} = 
    \begin{Bmatrix}
        p_0,p_1, ..., p_n = (x,y) \in G\\
        u_0,u_1,...,u_{n-1} \in \mathcal{U}\\
        t_0,t_1,...,t_{n-1} \in \mathbb{N} \setminus \{0,..,k-1\}
    \end{Bmatrix}, 
\end{equation*}
is a controlled chain between $p_0 \in E$ and $p_n = (x,y)$. Eventually, the chain $\mu_{(\delta,k)}$ reaches the region $A_1$ given by 
\begin{equation*}
    A_1 = \left\{(z,t) \in G: \frac{dm(x)+M(x)}{1-d^2} < t \leq \frac{dm(x)+M(x)}{1-d^2} +\frac{\delta+\varepsilon}{\cos{\alpha}}\right\}. 
\end{equation*}
since the last point of it is outside the $E$. 

As in the previous case, let us take $j_0 \in \{1,...,n-1\}$ such that $p_{j_0 + i} \notin E$ with $i \in \{0,...,n-j_0\}$ and $p_{j_0-1} \in E$. Taking $p_{j_0} = (p_{j_0}^1, p_{j_0}^2),$ we have
\begin{eqnarray}
    \sum_{j=0}^{t_{j_0} - 1} d^{t_{j_0}-1-j} p(u_j,p_{j_0}^1) + d^{t_0}p_{j_0}^2  \leq \frac{dm(x)+M(x)}{1-d^2} - d^{t_{j_0}} (f_{t_{j_0}}(p_{j_0}^1)) + d^{t_{j_0}}p_{j_0}^2.
\end{eqnarray}

By  (\ref{eq1}), (\ref{eq2}) and (\ref{eq3}), we have
\begin{eqnarray}
    f_k(p_{j_0}^1) = \left\{ 
    \begin{array}{ccc}
        \dfrac{dm(p_{j_0}^1)+M(p_{j_0}^1)}{1-d^2}, k \hbox{ if $k$ is even,}\\
        \dfrac{m(p_{j_0}^1)+dM(p_{j_0}^1)}{1-d^2},\hbox{ if $k$ is odd.}
    \end{array}
    \right.
\end{eqnarray}

As $\dfrac{dm(p_{j_0}^1)+M(p_{j_0}^1)}{1-d^2} < p_{j_0}^2 \leq \dfrac{dm(p_{j_0}^1)+M(p_{j_0}^1)}{1-d^2} + \dfrac{\varepsilon + \delta}{\cos{\alpha}}$, we have
\begin{equation*}
    S^{t_{j_0}}_{u_{j_0}}(p_{j_0}) \leq \frac{dm(p_{j_0}^1)+M(p_{j_0}^1)}{1-d^2} - d^{t_{j_0}} \left(f_{t_{j_0}}(p_{j_0}^1)\right) + d^{t_{j_0}}\left(\frac{dm(p_{j_0}^1)+M(p_{j_0}^1)}{1-d^2} + \frac{\varepsilon + \delta}{\cos{\alpha}}\right). 
\end{equation*}

If $t_{j_0}$ is even, then $ f_{t_{j_0}}(p_{j_0}^1) = \dfrac{dm(p_{j_0}^1)+M(p_{j_0}^1)}{1-d^2}$ and 
\begin{equation*}
    S^{t_{j_0}}_{u_{j_0}}(p_{j_0}) \leq \frac{dm(p_{j_0}^1)+M(p_{j_0}^1)}{1-d^2} + d^{t_{j_0}} \left(\frac{\varepsilon + \delta}{\cos{\alpha}}\right)\\
    <\frac{dm(p_{j_0}^1)+M(p_{j_0}^1)}{1-d^2}  + \frac{\varepsilon}{\cos{\alpha}}. 
\end{equation*}

If $t_{j_0}$ is odd, one can prove that 
$M(p_{j_0}^1) > m(p_{j_0}^1)$ implies  
\begin{equation*}
    - d^{t_{j_0}} \left( \frac{dM(p_{j_0}^1) + m(p_{j_0}^1)}{1-d^2} \right) + d^{t_{j_0}} \left( \frac{M(p_{j_0}^1) + dm(p_{j_0}^1)}{1-d^2} \right) < 0.
\end{equation*}

Then 
\begin{equation*}
    S^{t_{j_0}}_{u_{j_0}}(p_{j_0}) < \frac{dm(p_{j_0}^1)+M(p_{j_0}^1)}{1-d^2}  + \frac{\varepsilon}{\cos{\alpha}}. 
\end{equation*}

Following the same arguments we conclude that $p_{j_0 + i} \in A_1$, $n-1-j_0 \geq i \geq 0$ and, the $(\delta,k)-$controlled chain $\mu_{(\delta,k)}$ does not reach the open ball $B_{\delta}(x,y)$, a contradiction. $\blacksquare$

\begin{remark}The previous results  prove the existence of a unique chain control set for the restricted system $(\Sigma)$, when $h \equiv 1$ and $|d| < 1$.      
\end{remark}

\section{Controllability sets and chain controllability}

First we consider the case $d=1$, where the solution has the form 
\begin{equation*}
    \varphi(k,(x,y),u) = \left(x, \sum_{j=0}^{k-1}((a + g(u_j))x-a) + y\right). 
\end{equation*}

In this case the description of the control sets is summarized in the next proposition. 

\begin{proposition}\label{controlsetsd=1}If $d=1$, then for any $x \in \mathbb{R}_+$ satisfying 
\begin{equation}\label{inequalityinx}
    (a+g(u_m))x-a < 0 < (a+g(u_M))x-a, 
\end{equation}
the sets 
\begin{equation*}
    D_x=\{x\} \times \mathbb{R},
\end{equation*}
are control sets of $(\Sigma)$. 
\end{proposition}

\noindent\textit{Proof:} In fact, take $x \in \mathbb{R}_+$ and $y \in \mathbb{R}$. We have 
\begin{equation*}
    k[(a+g(u_m))x-a] + y  < y < k[(a+g(u_M))x-a] + y, k \in \mathbb{N}. 
\end{equation*}

By inequality (\ref{inequalityinx}), we get 
\begin{eqnarray*}
    k[(a+g(u_M))x-a] &\xrightarrow{}& \infty, \,\, \mbox{ for } k \rightarrow \infty \\
    k[(a+g(u_m))x-a] &\xrightarrow{}& -\infty, \,\, \mbox{ for } k \rightarrow \infty.
\end{eqnarray*}

Then 
\begin{eqnarray*}
    \mathcal{R}(x,y) &=& \{x\} \times \{[k((a+g(u_m))x-a) + y,k((a+g(u_M))x-a) + y]: k \in \mathbb{N}\}\\
    &=&\{x\} \times \mathbb{R} 
\end{eqnarray*}
for any pair $(x,y) \in G$ satisfying (\ref{inequalityinx}). Therefore, the set $\{x\} \times \mathbb{R}$ is a control set. $\blacksquare$

The following two results are purely technical, but necessary for our purposes. 

\begin{proposition}\label{maximalitychaind=1}If $d=1$ and $(a + g(u_M))x-a \leq 0$, then $\{x\} \times \mathbb{R}$ does not contain  control sets for any  $x \in \mathbb{R}_+$.
\end{proposition}

\noindent\textit{Proof:} Suppose that $D_x \subset \{x\} \times \mathbb{R}$ is a control set and $(x,y),(x,z) \in D_x$ with $y < z$, then 
\begin{equation*}
    \mathcal{R}(x,y) = \{x\} \times \left\{[k((a+g(u_m))x- a) + y,k((a+g(u_M))x-a) + y]: k \in \mathbb{N}\right\},
\end{equation*}
where $k((a+g(u_M))x-a) + y \leq y$. Given $\varepsilon = z - y$, we have
\begin{equation*}
    S^k_u(x,y) \leq y \iff z -S^k_u(x,y) \geq z - y = \varepsilon. 
\end{equation*}

Hence
\begin{equation*}
    ||(x,z) - (x,S^k_u(x,y))|| = |z - S^k_u(x,y)| \geq \varepsilon.
\end{equation*}

Therefore 
\begin{equation*}
    B_{\varepsilon}(x,z) \cap \mathcal{R}(x,y) = \emptyset,
\end{equation*}
which contradicts the approximate controllability of $D_x$. $\blacksquare$

\begin{figure}
    \centering
    \includegraphics[scale=0.7]{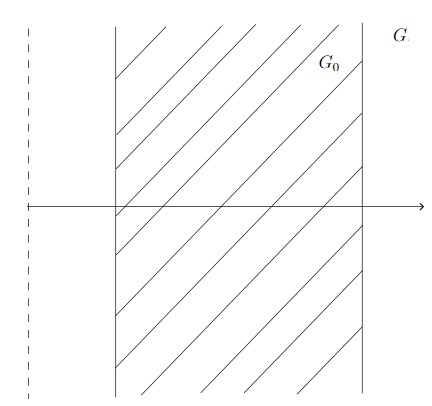}
    \caption{Region $G_0$ which contains every control set in the form $D_x = \{x\} \times \mathbb{R}$, with $d = 1$ and $(a+g(u_m))(a+g(u_M)) > 0$.}
    \label{g_0}
\end{figure}

\begin{proposition}If $d=1$ and $(a + g(u_m))x-a \geq 0$, then for any  $x \in \mathbb{R}_+$ the set $\{x\} \times \mathbb{R}$ does not contain  non-trivial control set. 
\end{proposition}

\noindent\textit{Proof:} If $D_x \subset \{x\} \times \mathbb{R}$ is a control set and $(x,y),(x,z) \in D_x$ with $y < z$, then 
\begin{equation*}
    \mathcal{R}(x,z) = \{x\} \times \left[k((a+g(u_m))x- a) + z,k((a+g(u_M))x-a) + z\right] 
\end{equation*}
where $k((a+g(u_M))x-a) + z \geq z$. Given $\varepsilon = z - y$, we have
\begin{equation*}
    S^k_u(x,z) \geq z > y.  
\end{equation*}

Then
\begin{equation*}
    ||(x,y) - (x,S^k_u(x,z))|| = |y - S^k_u(x,z)| = S^k_u(x,z)-y \geq z-y = \varepsilon.
\end{equation*}

Thus 
\begin{equation*}
    B_{\varepsilon}(x,y) \cap \mathcal{R}(x,z) = \emptyset,
\end{equation*}
which contradicts the approximate controllability of $D_x$. $\blacksquare$

Therefore, we can state the following result.   

\begin{proposition}\label{chaincontrolG_0}If $d=1$, the set 
\begin{equation*}
    G_0 = \{(x,y) \in \mathbb{R}_+ \times \mathbb{R}: (a+g(u_m))x-a \leq 0 \leq (a+g(u_M))x-a\}
\end{equation*}
satisfies the properties $1$ and $2$ of the definition (\ref{definchaincontrolset}). 
\end{proposition}

\noindent\textit{Proof:} We claim that the set of all $x_0$ satisfying (\ref{inequalityinx}) is bounded, if and only if, $(a+g(u_m))(a+g(u_M)) > 0$. In fact, consider the inequality
\begin{equation*}
    (a+g(u_m))x-a < 0 < (a+g(u_M))x-a.
\end{equation*}

If $a+g(u_m) = 0$, then $a+g(u_M) > 0$. Thus \begin{eqnarray*}
    \{x\in \mathbb{R}_+: (a+g(u_m))x-a < 0 < (a+g(u_M))x-a\} &=& \{x\in \mathbb{R}_+:-a < 0 < (a+g(u_M))x-a\}\\
    &=&  \left\{x\in \mathbb{R}_+: x > \frac{a}{a+g(u_M)}\right\}\\
    &=& \left(\frac{a}{a+g(u_M)},\infty\right). 
\end{eqnarray*} 

On the other hand, if $a+g(u_M) =0$, then $a+g(u_m) < 0$. Hence 
\begin{eqnarray*}
    \{x\in \mathbb{R}_+: (a+g(u_m))x-a < 0 < (a+g(u_M))x-a\} &=& \{x\in \mathbb{R}_+:(a+g(u_m))x-a < 0 <-a\}\\
    &=&  \left\{x\in \mathbb{R}_+: x > \frac{a}{a+g(u_m)}\right\}\\
    &=& \left(\frac{a}{a+g(u_m)},\infty\right). 
\end{eqnarray*}

Therefore, the set $\{x\in \mathbb{R}_+: (a+g(u_m))x-a < 0 < (a+g(u_M))x-a\}$ is unbounded. Reciprocally, if $(a+g(u_m))(a+g(u_M)) > 0$, then 
\begin{equation*}
    (a+g(u_m))x-a < 0 < (a+g(u_M))x-a,
\end{equation*}
implying that
\begin{eqnarray*}
    \frac{a}{a+g(u_M)} < &x& < \frac{a}{a+g(u_m)}, \hbox{ if }(a+g(u_m)) > 0\\
    \frac{a}{a+g(u_m)} < &x& < \frac{a}{a+g(u_M)}, \hbox{ if }(a+g(u_m)) < 0. 
\end{eqnarray*}

Therefore 
\begin{equation*}
    \{x \in \mathbb{R}_+: (a+g(u_m))x-a < 0 < (a+g(u_M))x-a\} 
\end{equation*}
is bounded. 

Now, we can prove the result. The property $1$ of the definition (\ref{definchaincontrolset}) is trivial, since the set $G_0$ contains every line $\{x\} \times \mathbb{R}$ satisfying $(a+g(u_m))x-a \leq 0 \leq (a+g(u_M))x-a$. For the chain controllability, take $(x,y), (z,t) \in \hbox{int}(G_0)$ and $(\varepsilon,k) \in \mathbb{R}_+ \times \mathbb{N}$. If $z\in (x,x+\varepsilon)$, consider $(z,t_1) \in G_0$ such that $\varphi(k,(x,y),u) \in B_{\varepsilon}(z,t_1)$ for some arbitrary $u \in \mathcal{U}$. As $\{x\} \times \mathbb{R}$ is a control set, there are $u_1 \in \mathcal{U}$ and $T_1 \geq k$ such that 
\begin{equation*}
    \varphi(T_1, (z,t_1),u_1) \in B_{\varepsilon}(z,t). 
\end{equation*}

Therefore, the set 
\begin{equation*}
    \eta_{(\varepsilon,k)} = 
    \begin{Bmatrix}
        (x,y), (z,t_1), (z,t) \in G\\
        u,u_1 \in \mathcal{U}\\
        k, T_1 \in \mathbb{N}_0 \setminus \{0,...,k-1\}
    \end{Bmatrix}
\end{equation*}
is a controlled chain between $(x,y)$ and $(z,t)$. If $z \geq x+\varepsilon$, take $x_1 \in (x,x+\varepsilon)$. Then we can construct a $(\varepsilon,k)-$controlled chain between $(x,y)$ and $(x_1,y)$. If $z \in (x_1, x_1+  \varepsilon)$, one can apply the same initial argument to get a $(\varepsilon,k)-$controlled chain between $(x_1,y)$ and $(z,t)$. Then, concatenating such chains, there is a $(\varepsilon,k)-$controlled chain between $(x,y)$ and $(z,t)$. Applying this successively, one can construct a controlled chain between $(x.y)$ and $(z,t)$.

Now, considering $(x,y) \in \partial G_0$ and $G_0$ bounded, we have by initial claim that $x = \frac{a}{a+g(u_m)}$ or $x = \frac{a}{a+g(u_M)}$. If $x = \frac{a}{a+g(u_m)}$, the reachable set is  
\begin{equation*}
   \mathcal{R}(x,y) =  \{x\} \times \left\{[ y,k((a+g(u_M))x-a) + y]: k \in \mathbb{N}\right\}.  
\end{equation*}

Fix $(\varepsilon,k) \in \mathbb{R}_+ \times \mathbb{N}$ and consider $(x_1,y_1) \in D_{x_1} \subset \hbox{int}G_0$ such that $\varphi(k,(x,y),u) \in B_{\varepsilon}(x_1,y_1)$, then we can connect by a $(\varepsilon,k)-$controlled chain any two points of $G_0$. The same for $x = \frac{a}{a+g(u_M)}$. Now, if $G_0$ is unbounded, then $a+g(u_m) = 0$ or $a+g(u_M) = 0$. For the case $a+g(u_m) = 0$, we have $G_0 = \{(x,y) \in G: x \geq \frac{a}{a+g(u_M)}\}$. In the same way we prove that $G_0$ is chain controllable. The argument is similar for $a+g(u_M) = 0$. Then $G_0$ is chain controllable.$\blacksquare$

Considering the case $d = -1$, one can prove by induction that 
\begin{equation*}
    \sum_{j=0}^{k-1} (-1)^{k-1-j} p(u_j,x) = \left\{
    \begin{array}{cc}
        \sum_{j=0}^{\frac{k}{2}-1} (-1)^{k-1-2j} p(u_{2j},x) + \sum_{j=0}^{\frac{k}{2}-1} (-1)^{k-2-2j} p(u_{2j+1},x), \hbox{ $k\geq 2$ and even},\\
        \sum_{j=0}^{\frac{k-1}{2}} (-1)^{k-1-2j} p(u_{2j},x) + \sum_{j=0}^{\frac{k-1}{2}-1} (-1)^{k-2-2j} p(u_{2j+1},x), \hbox{ $k\geq 3$ and odd}.
    \end{array}
    \right.
\end{equation*}

Then 
\begin{equation*}
\begin{array}{ccl}
    \frac{k}{2}(m(x)-M(x))\leq \sum_{j=0}^{k-1} (-1)^{k-1-j} p(u_j,x) \leq \frac{k}{2}(M(x)-m(x)),& k \hbox{ even },\\
    (\frac{k-1}{2})(m(x)-M(x)) + m(x)\leq \sum_{j=0}^{k-1} (-1)^{k-1-j} p(u_j,x) \leq (\frac{k-1}{2})(M(x)-m(x)) + M(x),& k \hbox{ odd }.
\end{array}
\end{equation*}

Finally, we get the following result that ensures the chain controllability of $(\Sigma)$ when $d=-1$. 

\begin{proposition}If $d=-1$, the subsets 
\begin{equation*}
    D_x=\{x\}\times \mathbb{R}, \mbox{ with } x \in \mathbb{R}_+,
\end{equation*}
are control sets and $G$ is chain controllable. 
\end{proposition}

\noindent\textit{Proof:} Take $x \in \mathbb{R}_+$ and recall that $M(x) > m(x)$ for every $x \in \mathbb{R}_+$. Considering the sequences
\begin{equation*}
    x_k =\left\{
    \begin{array}{ll}
        \frac{k}{2}(m(x)-M(x)),& k\hbox{ even }\\
        (\frac{k-1}{2})(m(x)-M(x)),& k \hbox{ odd}
    \end{array}
    \right.
\end{equation*}
and
\begin{equation*}
    y_k =\left\{
    \begin{array}{ll}
        \frac{k}{2}(M(x)-m(x)),& k\hbox{ even }\\
        (\frac{k-1}{2})(M(x)-m(x)),& k \hbox{ odd}
    \end{array}
    \right.
\end{equation*}
for $k \geq 2$, it is not hard to prove that $x_k \longrightarrow -\infty$ and $y_k \longrightarrow \infty$. Given $(x,y) \in G$, we have
\begin{equation*}
    \varphi(k,(x,y),u) = (x, \sum_{j=0}^{k-1}(-1)^{k-1-j} p(u_j,x) + (-1)^{k}y) 
\end{equation*}
with $\sum_{j=0}^{k-1}(-1)^{k-1-j} p(u_j,x)$ satisfying 
\begin{equation*}
\begin{array}{rccccl}
    \frac{k}{2}(m(x)-M(x))&\leq& \sum_{j=0}^{k-1} (-1)^{k-1-j} p(u_j,x) &\leq& \frac{k}{2}(M(x)-m(x)),& k \hbox{ even, }\\
    (\frac{k-1}{2})(m(x)-M(x)) + m(x)&\leq& \sum_{j=0}^{k-1} (-1)^{k-1-j} p(u_j,x) &\leq &(\frac{k-1}{2})(M(x)-m(x)) + M(x),& k \hbox{ odd }.
\end{array}
\end{equation*}

Now, fix $x \in \mathbb{R}_+$. For $z,y \in \mathbb{R}_+$, there is a $k \in \mathbb{N}$ even, such that 
\begin{equation*}
    z \in \left[ \frac{k}{2}(m(x)-M(x)) + y,\frac{k}{2}(M(x)-m(x)) + y \right]. 
\end{equation*}

Besides, there is a $t \in \mathbb{N}$ odd, such that 
\begin{equation*}
    z \in \left[ \frac{t-1}{2}(m(x)-M(x)) + m(x)- y,\frac{t-1}{2}(M(x)-m(x)) + M(x) - y \right].
\end{equation*}

In both cases, $(x,z) \in \mathcal{R}(x,y)$. The maximality and invariance are trivial. Then the chain controllability follows by the same argument as in Proposition \ref{chaincontrolG_0}. $\blacksquare$

\textbf{Acknowledgement. }We appreciate the careful reading and the
constructive comments of an anonymous reviewer which helped to improve the paper.

\section{Declarations}

\subsection{Funding}

No funding received

\subsection{Conflict of interest/Competing interests}

No interests of a financial or personal nature

\subsection{Ethics approval}

Not applicable

\subsection{Consent to participate}

According to \newline
https://www.springer.com/gp/editorial-policies/informed-consent this item it
is not applicable

\subsection{Consent for publication}

According to \newline
https://www.springer.com/gp/editorial-policies/informed-consent this item it
is not applicable

\subsection{Availability of data and materials}

Not applicable

\subsection{Code availability}

Not applicable

\subsection{Authors' contributions}

All authors contributed to all sections. All authors reviewed the final manuscript.


\begin{thebibliography}{9999}

\bibitem{sontalb}Albertini, F., Sontag, E.\emph{ Discrete-time transitivity and accessibility: analytic systems } SIAM J. Control and Optimization. 1993. Vol. 31, No. 6, pp. 1599-1622.

\bibitem{sontalb1} Albertini, F., Sontag, E.\emph{ Some connectons between chaotic dynamical systems and control systems } Proc. European Control Conf. 1991. Vol 1. Grenoble. 58-163, 

\bibitem{AyalaeAdriano}Ayala, V., Da Silva, A. \emph{ The control set of a linear control system on the two dimensional solvable Lie group. } Journal of Differential Equations. 2020.

\bibitem{TAJ}Cavalheiro, T. ; Santana, A. J.; Cossich, J. A. N. \emph{Controllability of discrete-time linear
systems on solvable Lie groups. } arXiv preprint arXiv:2302.00145, (2023).

\bibitem{CCS}Colonius, F., Cossich, J. A. N., Santana, A. J. \emph{ Controllability Properties and Invariance Pressure for Linear
Discrete-Time Systems. }J Dyn Diff Equat. 2022. 34, 5–28.

\bibitem{CCS1}Colonius, F., Cossich, J. A. N., Santana, A. J. \emph{ Outer invariance entropy for discrete-time linear
systems on Lie groups. } ESAIM: Control, Optimisation and Calculus of Variations. 2021. vol. 27.

\bibitem{CCS2}Colonius, F., Kliemann, W. \emph{ The dynamics of control. } Springer. New York. (1999).


\bibitem{sontag1}Sontag, E.D., \emph{ Mathematical Control Theory. Deterministic finite - dimensional systems, } Springer-Verlag, New York, 1998.

\end{thebibliography}
\end{document}